\documentclass[12pt, a4paper]{article}

\usepackage{geometry}
\usepackage{float}
\usepackage{graphicx}
\usepackage{color}
\usepackage[pagewise]{lineno}
\usepackage{ulem}
\usepackage{multirow}
\usepackage{booktabs}
\usepackage{authblk}
\usepackage{amsmath, amssymb}
\usepackage[pagewise]{lineno}
\usepackage{indentfirst}
\usepackage[marginal]{footmisc}
\usepackage{ntheorem}

\theorembodyfont{\normalfont}  

\newtheorem{remark}{Remark}

\makeatletter 
\@addtoreset{equation}{section} 
\makeatother
\usepackage[colorlinks,
            linkcolor=red,
            anchorcolor=red,
            citecolor=red
            ]{hyperref}
\DeclareMathOperator*{\argmin}{argmin}

\newcommand{\mcd}{\mathcal{D}}

\newcommand{\rd}{\mathrm{d}}

\newcommand{\bn}{\boldsymbol{n}}
\newcommand{\boeta}{\boldsymbol{\eta}}
\newcommand{\bx}{\boldsymbol{x}}
\newcommand{\bxi}{\boldsymbol{\xi}}
\newcommand{\abs}[1]{\left\vert#1\right\vert}
\newcommand{\brac}[1]{\left(#1\right)}

\newcommand{\norm}[1]{\left\Vert#1\right\Vert}

\newcommand{\nn}{\nonumber}

\title{The Quadratic Wasserstein Metric With Squaring Scaling For Seismic Velocity Inversion}
\author[1]{Zhengyang Li}
\author[2]{Yijia Tang}
\author[$*$3]{Jing Chen}
\author[$*$1]{Hao Wu}
\affil[1]{Department of Mathematical Sciences, Tsinghua University, Beijing, China 100084.}
\affil[2]{School of Mathematical Sciences, Shanghai Jiao Tong University, Shanghai, China 200240.}
\affil[3]{Division of mathematical sciences, school of physical and mathematical sicences, Nanyang Technological University, Singapore 639798.}

\begin{document}
\maketitle
\footnote{\boldmath{$^*$} \textbf{Corresponding author}. \\}

\begin{abstract}
The quadratic Wasserstein metric has shown its power in measuring the difference between probability densities, which benefits optimization objective function with better convexity and is insensitive to data noise. Nevertheless, it is always an important question to make the seismic signals suitable for comparison using the quadratic Wasserstein metric. The squaring scaling is worth exploring since it guarantees the convexity caused by data shift. However, as mentioned in [Commun. Inf. Syst., 2019, 19:95-145], the squaring scaling may lose uniqueness and result in more local minima to the misfit function. In our previous work [J. Comput. Phys., 2018, 373:188-209], the quadratic Wasserstein metric with squaring scaling was successfully applied to the earthquake location problem. But it only discussed the inverse problem with few degrees of freedom. In this work, we will present a more in-depth study on the combination of squaring scaling technique and the quadratic Wasserstein metric. By discarding some inapplicable data, picking seismic phases, and developing a new normalization method, we successfully invert the seismic velocity structure based on the squaring scaling technique and the quadratic Wasserstein metric. The numerical experiments suggest that this newly proposed method is an efficient approach to obtain more accurate inversion results.


\end{abstract}

\section{Introduction}
\label{sec:intro}
Full waveform inversion (FWI) has been receiving wide attention in recent years \cite{DuYa:20,EnYa:21, MeBrViOp:13, WaYaJiWu:19, YaEnSuFr:18, YaEn:18} due to its high-resolution imaging in geophysical properties. Generally, it can be formulated as a PDE constrained optimization problem in mathematics, which consists of two parts \cite{ViOp:09}: the forward modeling of seismic wavefield, and the optimization problem searching for suitable model parameters to minimize the mismatch between the predicted and observed seismic signals. In previous decades, limited by the computing power, most tomography methods were based on the ray theory, which ignores finite frequency phenomena such as wave-front healing and scattering \cite{HuDaNo:01}, and thus results in low-resolution inverion results. With the rapid development of computing power and the forward modeling method, more accurate synthetic signals could be computed by directly simulating seismic wave propagation. This makes it possible to obtain high-resolution results by FWI, which could provide important information for seismic hazard assessment \cite{TaLiuMaTr:10} and exploration geophysics \cite{ViOp:09}.

The $L^2$ metric-based model is the simplest and most common FWI. However, it suffers from the well-known cycle skipping problem \cite{ViOp:09} that the solution may be trapped in the local minima during the iteration, leading to incorrect inversion results. The quadratic Wasserstein metric ($W_2$) from the Optimal transport (OT) theory \cite{Vi:03, Vi:08} seems to be a solution to the above problem. It measures the difference between two probability distributions by minimizing the transport cost from one distribution to the other, which is insensitive to the data noise and keeps convexity to the data shift, dilation, and partial amplitude change \cite{EnFr:14, EnFrYa:16}. The number of local minima of the FWI model based on this metric is therefore significantly reduced. Thus, it is favored by researchers and has been widely applied to the earthquake location and seismic tomography \cite{ChChWuYa:18, EnFr:14, EnFrYa:16, EnYa:19a, EnYa:19b, YaEnSuFr:18, YaEn:18, ZhChWuYaQi:21}. In applying the quadratic Wasserstein metric to the seismic inverse problem, there is a critical problem. The quadratic Wasserstein metric compares the normalized and nonnegative data while the seismic signal does not meet this requirement. Thus, various techniques are developed to deal with this problem, e.g., linear scaling \cite{YaEnSuFr:18}, squaring scaling \cite{ChChWuYa:18}, and exponential scaling \cite{QiJaVaYaEn:17}. Among all these methods, squaring scaling is considered to maintain the convexity of the optimization objective function. But this method seems to lose uniqueness and result in additional minima. This may be the reason why we haven't seen the application of squaring scaling and quadratic Wasserstein metric to the velocity inversion problem. Moreover, there are also some other metrics based on the OT theory, e.g., the WFR metric and the KR norm, which have been successfully applied to the seismic inverse problem \cite{MeBrMeOuVi:16a, MeBrMeOuVi:16b, ZhChWuYaQi:21}.

In our previous work \cite{ChChWuYa:18}, the quadratic Wasserstein metric with squaring scaling is successfully applied to the earthquake location problem. The squaring scaling ensures the differentiability and nice convexity property, leading to a large convergent domain and accurate inversion results. However, it is still a challenging problem for velocity inversion with a large number of degrees of freedom since the squaring scaling may lose uniqueness and result in additional local minima to the misfit function \cite{EnYa:19a}. In this work, we would like to provide a comprehensive approach to the seismic velocity inversion based on squaring scaling and the quadratic Wasserstein metric. The key ingredient of this work consists of two parts. First, for seismic velocity inversion, the fundamental geophysical characteristic of seismic signals should be taken into account. For example, certain erroneous seismic signals and multi-arrival seismic signals, which have destructive effects on the inverse process, should be deleted in the preprocessing stage. Moreover, a more accurate optimal transport map can be obtained by picking appropriate seismic phases. Secondly, a new normalization method is developed to obtain a more accurate optimal transport map for the squared seismic signals. From this, we can calculate better sensitivity kernels, which are more consistent with physical intuition.

The rest of the paper is organized as follows. In Section \ref{sec:overview}, we briefly review the mathematical formula of seismic velocity inversion and the basics of the quadratic Wasserstein metric. We discuss important issues in the inversion and present detailed implementations in Section \ref{sec:proc}. Meanwhile, we illustrate the necessity of our method by some toy models. In Section \ref{sec:num}, the numerical experiments are provided to demonstrate the effectiveness and efficiency of our method. Finally, we conclude the paper in Section \ref{sec:con}.

\section{The quadratic Wasserstein metric and seismic velocity inversion}
\label{sec:overview}

We review the full waveform seismic tomography and the adjoint state method in this section. The mathematical formulation of seismic velocity inversion can be written as the PDE constrained optimization problem,
\begin{equation} \label{eqn:inv_prob}
	c_T(\bx)=\argmin_{c(\bx)} \Xi(c(\bx)),\quad \Xi(c(\bx))=\sum_{i=1}^{N}\sum_{j=1}^{M}\chi_{ij}(c(\bx)),
\end{equation}
where index $(i,j)$ indicates the source-receiver pair. We used $N$ seismic events, and considered $M$ seismic signals for each event. Correspondingly, the misfit function $\chi_{ij}$ is defined as
\begin{equation} \label{eqn:mist_fun}
	\chi_{ij}(c(\bx))=\mcd(s_{ij}(t;c(\bx)),d_{ij}(t)).
\end{equation}
Here, $\mcd$ is the distance function that measures the difference between the real seismic signal $d_{ij}(t)$ and the synthetic signal $s_{ij}(t;c(\bx))$, which can be regarded as the solution
\begin{equation} \label{eqn:real_syn}
	d_{ij}(t)=u_i(\boeta_j,t;c_T(\bx)), \quad s_{ij}(t;c(\bx))=u_i(\boeta_j,t;c(\bx)),
\end{equation}
of the following acoustic wave equation with the initial boundary condition

\begin{align}  
	& \frac{\partial^2 u_i(\bx,t;c(\bx))}{\partial t^2}=\nabla\cdot\brac{c^2(\bx)\nabla u_i(\bx,t;c(\bx))}
		+R(t-\tau_i)\delta(\bx-\bxi_i), \quad \bx\in\Omega, t>0, \label{eqn:wave} \\
	& u_i(\bx,0;c(\bx))=\frac{\partial u_i(\bx,0;c(\bx))}{\partial t}=0, \quad \bx\in\Omega, \label{ic:wave} \\
	& \bn\cdot \brac{c^2(\bx)\nabla u_i(\bx,t;c(\bx))}=0, \quad \bx\in\partial \Omega, t>0. \label{bc:wave}
\end{align} 
Here, the locations of the earthquake and receiver station are $\bxi_i$ and $\boeta_j$, the origin time of the earthquake is $\tau_i$. The seismic rupture is modeled by the point source $\delta(\bx-\bxi)$ since its scale is much smaller compared to the scale of seismic wave propagation \cite{AkRi:80,Ma:15}. And the source time function is simplified as the Ricker wavelet 
\begin{equation} \label{eqn:ricker}
R(t)=A \left( 1-2 \pi^2f_0^2t^2 \right) e^{- \pi^2f_0^2t^2},
\end{equation}
where $f_0$ denotes the dominant frequency, and $A$ is the normalization factor. The outward unit normal vector to the simulation domain boundary $\partial \Omega$ is $\bn$. In practice, the perfectly matched layer absorbing boundary condition \cite{KoTr:03} is used to deal with the propagation of waves outside the area. In this section, we use the reflection boundary condition to simplify the derivation.

\begin{remark}
Here, we consider the trace by trace strategy \cite{YaEnSuFr:18} to apply the 1-D quadratic Wasserstein metric to the waveform inversion. Considering the fact that receiver stations are located far from each other on the geological scale, this approach is more in line with physical reality and also easier in mathematics.
\end{remark}

\subsection{The adjoint method}  \label{sec:adjoint}
Below, we briefly review the adjoint method \cite{EnFrYa:16, Pl:06} for solving the optimization problems \eqref{eqn:inv_prob}-\eqref{eqn:ricker}. For small perturbation of seismic velocity structure $\delta c$, it causes the perturbation of the wavefield

\begin{equation} \label{rel:pert_s}
	\delta u_i(\bx,t;c(\bx))=u_i(\bx,t;c+\delta c)-u_i(\bx,t;c).
\end{equation}
For the sake of brevity, we will omit the parameter $c(\bx)$ of the wavefield and the signals in the following. The perturbation $\delta u_i(\bx,t)$ satisfies the equations
 
\begin{align} 
	& \frac{\partial^2 \delta u_i(\bx,t)}{\partial t^2}=\nabla\cdot\brac{c^2(\bx)\nabla \delta u_i(\bx,t)}  \label{eqn:pert_s}\\
	&\quad \quad \quad \quad \quad\ \ +\nabla\cdot\brac{\brac{2c(\bx)+\delta c(\bx)}\delta c(\bx)\nabla (u_i+\delta u_i)(\bx,t)}, \quad \bx\in\Omega, \nn \\
	& \delta u_i(\bx,0)=\frac{\partial \delta u_i(\bx,0)}{\partial t}=0, \quad \bx\in\Omega, \label{ic:pert_s} \\
	& \bn\cdot \brac{c^2(\bx)\nabla \delta u_i(\bx,t)+\brac{2c(\bx)+\delta c(\bx)}\delta c(\bx)\nabla (u_i+\delta u_i)(\bx,t)}=0, \quad \bx\in\partial \Omega. \label{bc:pert_s}
\end{align}
Multiply test function $w_i(\bx,t)$ on equation \eqref{eqn:pert_s} and integrate it on $\Omega\times[0,t_f]$ for sufficient large time $t_f$. Using integration by parts yields

{\small\begin{multline} \label{eqn:integration}
	\int_0^{t_f}\int_{\Omega}\frac{\partial^2 w_i}{\partial t^2}\delta u_i\rd \bx\rd t
		-\int_{\Omega}\left.\frac{\partial w_i}{\partial t}\delta u_i\right|_{t=t_f}\rd \bx
		+\int_{\Omega}\left.w_i\frac{\partial \delta u_i}{\partial t}\right|_{t=t_f}\rd \bx \\
	=\int_0^{t_f}\int_{\Omega}\nabla\cdot(c^2\nabla w_i)\delta u_i\rd \bx \rd t
		-\int_0^{t_f}\int_{\partial \Omega}\bn\cdot (c^2\nabla w_i)\delta u_i\rd \zeta \rd t
        -\int_0^{t_f}\int_{\Omega}\brac{2c+\delta c}\delta c\nabla w_i\cdot\nabla (u_i+\delta u_i)\rd \bx\rd t \\
	\approx\int_0^{t_f} \int_{\Omega}\nabla\cdot(c^2\nabla w_i)\delta u_i\rd \bx \rd t
         -\int_0^{t_f} \int_{\partial \Omega}\bn\cdot (c^2\nabla w_i)\delta u_i\rd \zeta \rd t
         -\int_0^{t_f}\int_{\Omega}2c\delta c\nabla w_i\cdot\nabla u_i\rd \bx\rd t,
\end{multline}}
where the higher-order terms are ignored in the last step since we can naturally assume that $\norm{\delta u_i}\ll \norm{u_i}$ and $\norm{\delta c(\bx)}\ll \norm{c(\bx)}$.

On the one hand, the perturbation of misfit $\delta \chi_{ij}$ results from the wave speed perturbation $\delta c(\bx)$, which writes

\begin{multline*}
	\delta \chi_{ij}(c)= \mcd\big(s_{ij}(t)+\delta s_{ij}(t),d_{ij}(t)\big)-\mcd\big(s_{ij}(t),d_{ij}(t)\big) \\
	\approx \langle Q_{ij}(t), \; \delta s_{ij}(t) \rangle
	=\int_0^{t_f}Q_{ij}(t) \delta s_{ij}(t)\rd t.
\end{multline*}
Here, $Q_{ij}(t)$ indicates the Fr\'echet gradient of the distance $\mcd$ with respect to the synthetic data $s_{ij}(t)$:

\begin{equation}\label{equ:adjoint_source}
Q_{ij}(t)=\nabla_s \mcd(s,d)\big|_{s=s_{ij}(t),d=d_{ij}(t)},
\end{equation}
which will be specified later. Let $w_i(\bx,t)$ satisfy the adjoint equation

\begin{align} 
	& \frac{\partial^2 w_i(\bx,t)}{\partial t^2}=\nabla\cdot\brac{c^2(\bx)\nabla w_i(\bx,t)} 
		+\sum_{j=1}^{M} Q_{ij}(t)\delta (\bx-\boeta_j), \quad \bx\in\Omega, \label{eqn:adjoint} \\
	& w_i(\bx,t_f)=\frac{\partial w_i(\bx,t_f)}{\partial t}=0, \quad \bx\in\Omega, \label{ic:adjoint} \\
	& \bn\cdot\brac{c^2(\bx)\nabla w_i(\bx,t)}=0, \quad \bx\in\partial \Omega. \label{bc:adjoint}
\end{align}
Multiply $\delta u_i(\bx,t)$ on equation \eqref{eqn:adjoint}, integrate it on $\Omega\times[0,t_f]$ and subtract \eqref{eqn:integration} to obtain

\begin{multline*}
	\sum_{j=1}^{M}\int_0^{t_f}Q_{ij}(t) \delta s_{ij}(t) \rd t = 
	\sum_{j=1}^{M}\int_0^{t_f} \int_\Omega Q_{ij}(t)  \delta(\bx-\boeta_j) \delta u_i(\bx,t) \rd t \\
	=-\int_0^{t_f}\int_{\Omega}2c(\bx)\delta c(\bx)\nabla w_i(\bx,t)\cdot\nabla u_i(\bx,t)\rd \bx\rd t.
\end{multline*}
The linear relationship between $\delta \Xi$ and $\delta c(\bx)$ is established as

\begin{equation} \label{eqn:delta}
	\delta \Xi(c)=\sum_{i=1}^{N}\sum_{j=1}^{M}\delta \chi_{ij}(c)=\sum_{i=1}^{N}\int_{\Omega}K_i(\bx)\delta c(\bx)\rd\bx,
\end{equation}
where the sensitivity kernel of the $i$-th source for $c(\bx)$ is defined as
\begin{equation} \label{eqn:kernel}
	K_i(\bx)=-\int_0^{t_f}2c(\bx)\nabla w_i(\bx,t)\cdot\nabla u_i(\bx,t)\rd t.
\end{equation}

\subsection{The quadratic Wasserstein metric} \label{sec:Wasserstein}
As we discussed at the beginning of this section, the synthetic signal $s_{ij}(t)$ and real seismic signal $d_{ij}(t)$ are time series. As we know, the quadratic Wasserstein metric between the 1-D probability density functions has an analytic form \cite{ChChWuYa:18, Vi:03, Vi:08, YaEnSuFr:18}, i.e.,

\begin{equation} \label{eqn:qwm_sol}
W_2^2(f,g)=\int_0^{t_f}\abs{t-T(t)}^2 f(t)\rd t,\quad T(t)=G^{-1}\brac{F(t)}.
\end{equation}
Here $f(t),\; g(t)$ are probability density functions defined on $[0,t_f]$ and $F(t),\; G(t)$ are cumulative density functions defined on $[0,t_f]$,

\begin{equation*}
F(t)=\int_{0}^{t_f} f(\tau)\rd\tau,\quad G(t)=\int_{0}^{t_f} g(\tau)\rd\tau.
\end{equation*}

Note that the seismic signals are not probability density functions. We need to transform them into nonnegative and normalized functions for the quadratic Wasserstein metric comparison. In other words, the misfit function defined in \eqref{eqn:mist_fun} can be written as

\begin{equation} \label{eqn:ref_dis}
\chi_{ij}=\mcd(s_{ij}(t),d_{ij}(t))=W^2_2(\mathcal{P}(s_{ij}(t)),\mathcal{P}(d_{ij}(t))).
\end{equation}
The operator $\mathcal{P}$ converts the seismic signals into probability density functions, including processing them into nonnegative and normalized time series. In the later part, we will discuss this in detail. Thus, we can obtain the expression of the Fr\'echet gradient \cite{ChChWuYa:18, YaEnSuFr:18} mentioned in \eqref{equ:adjoint_source},

\begin{equation} \label{equ:frechet_gradient}
	\nabla_s \mcd(s,d)=\nabla_f W^2_2(f,g)|_{f=\mathcal{P}(s),g=\mathcal{P}(d)} \cdot \nabla_s \mathcal{P}(s)
	=\left\langle 2\int_0^t \tau-T(\tau) \rd \tau ,\nabla_s \mathcal{P}(s) \right\rangle.
\end{equation}

\section{Data preprocessing and new normalization}\label{sec:proc}
In this section, we discuss two important issues when carrying out seismic velocity inversion. First of all, when using real data for inversion, we do not use all the data in each iteration. Some data, such as the case where the direct wave and the reflected wave arrive simultaneously, are difficult to use and can be ignored. In order to avoid the mismatch between different types of seismic phases, we only retain the direct waves in the real seismic signals and the synthetic signals. This processing procedure ensures reasonable optimal transport maps and accurate sensitivity kernels. Secondly, we will carefully design the operator $\mathcal{P}$ to get a better OT map $T$. In the following, we will present detailed implementations and discussions.

\subsection{Selecting source-receiver pairs and picking seismic phases} 
The complex subsurface structures, such as the velocity discontinuity interfaces, may lead to different types of seismic phases, including the direct wave and the reflected wave. These seismic waves propagate along different wave paths and carry distinct underground structure information. Sometimes, the direct wave and the reflected wave arrive simultaneously and can not be distinguished, called the multipath phenomenon \cite{RaSaHa:10}. It is not trivial to extract robust information from this kind of constraint. In practice, these source-receiver pairs are always manually excluded to avoid interference caused by unreliable constraints \cite{ChBa:06, HuYaTo:16}. We will also use this strategy in this study.

From the perspective of signal processing, different phases of the real seismic signal and the synthetic signal should be matched separately. If there is a matching error, for example, part of the direct wave of the synthetic signal is matched with part of the reflected wave of the real seismic signal, it would lead to the optimal transport map being inconsistent with basic seismic knowledge and further result in the artifacts in the sensitivity kernel \cite{EnFr:14}. In particular, for the squaring scaling and quadratic Wasserstein metric based seismic velocity inversion, this problem is more prominent. The reason is that the quadratic Wasserstein metric requires mass conservation and global match. When the masses of the real seismic signal and the synthetic signal are unbalanced in the same phase, the mass transportation between different phases will occur, causing the inconsistency between the OT map with seismic reality. Moreover, the squaring scaling could further magnify the problem. The idea of solving the above problems is also easy. By picking the phases, we only match the same phases of the real seismic signals and the synthetic signals. This is a common strategy in seismic inversion \cite{MaTaChChTr:09,ChHiLe:17}, and it can be achieved simply by calculating the arrival time of the direct phase and the reflected phase \cite{ChNiPi:14,WaTk:20}.

\begin{figure}[H]
	\centering
	\includegraphics[width=1\textwidth]{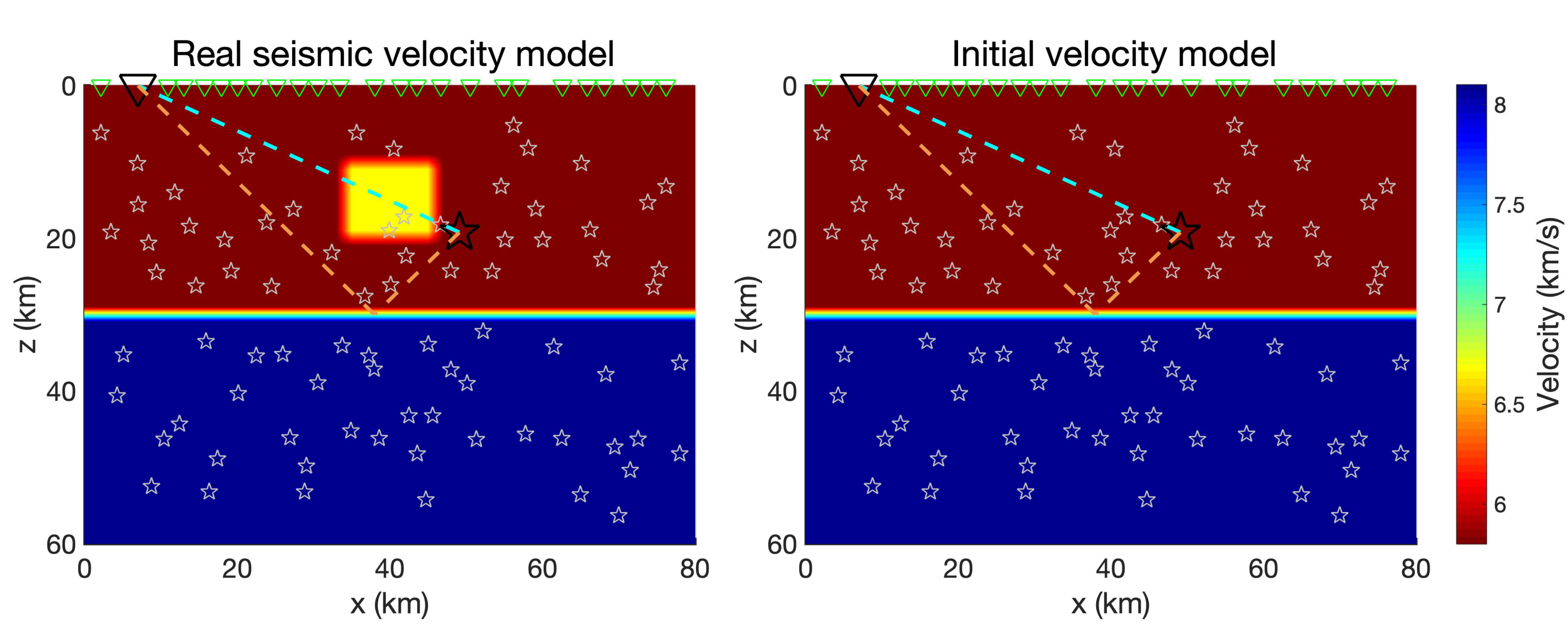}
	\caption{Illustration of the two-layer model. Left: the real seismic velocity model with a high-velocity anomaly; Right: the initial velocity model. The green inverted triangles indicate the receiver stations and the white stars indicate the earthquakes. The specific source-receiver pair is highlighted by the black star and inverted triangle. The cyan and tan dashed lines are the direct wave path and the reflected wave path, respectively.} 
	\label{fig:section3_velocity_52_2}
\end{figure}

Next, we explain the necessity of the above-mentioned data preprocessing method. The initial and real seismic velocity models are shown in Figure \ref{fig:section3_velocity_52_2}, and the parameter settings can be found in Section 4.1. The main goal is to detect the high-velocity anomaly above the Moho discontinuity.

Whether initial or real seismic velocity models, there are at least two paths from the earthquake hypocenter to the receiver station: the direct wave (cyan dashed lines) and the reflected wave (tan dashed lines). In the real seismic velocity model, the wave amplitude of the direct wave signal is slightly smaller since it partially reflects when passing through the high-velocity anomaly. On the other hand, the reflected wave signal should be the same since the velocity structure on the reflected wave path is the same in the initial and real seismic velocity models, see Figure \ref{fig:section3_velocity_52_2} for illustration.

\begin{figure}[H]
	\centering
	\includegraphics[width=1\textwidth]{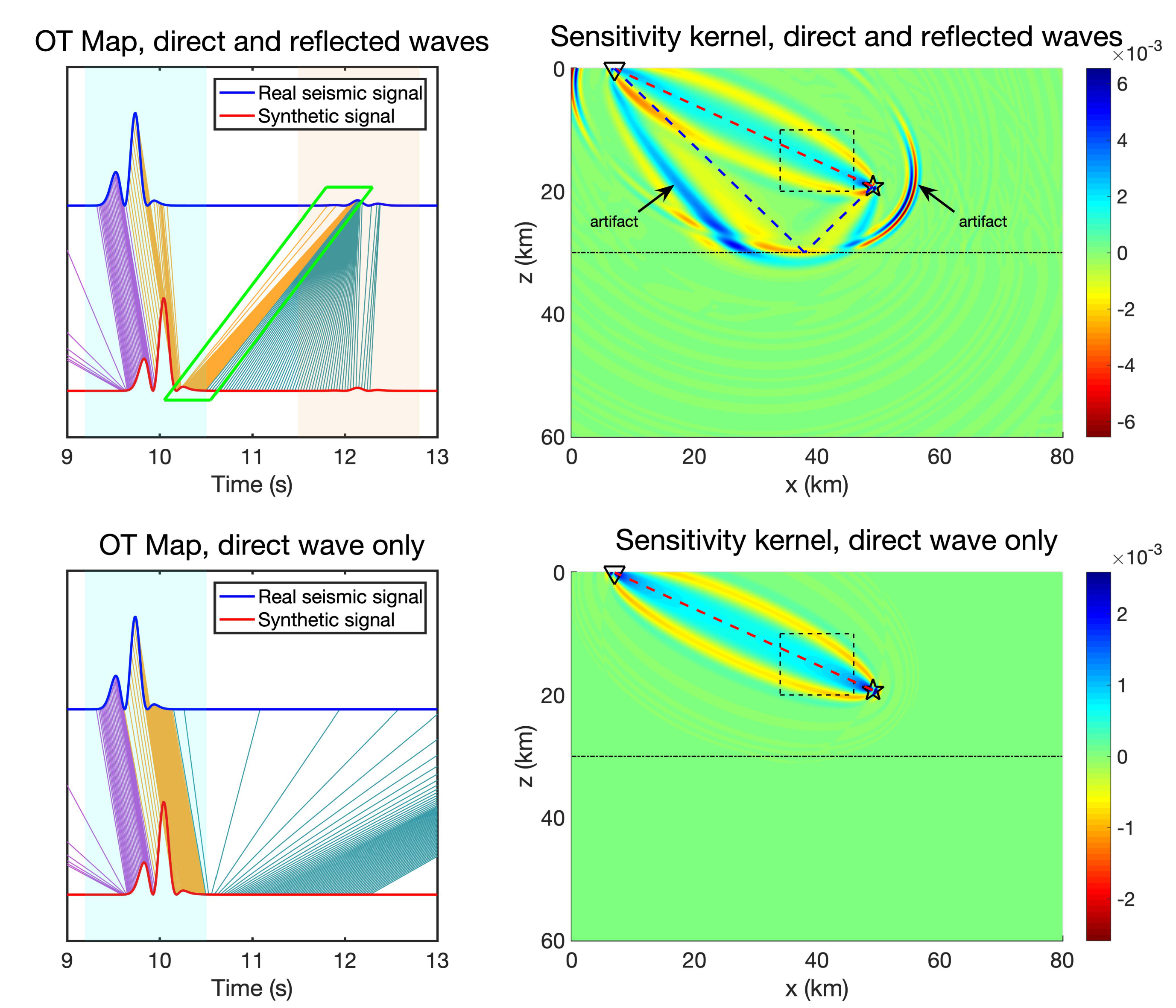}
	\caption{Illustration of the Optimal Transport map between the real seismic signal and synthetic signal (left) and the sensitivity kernel (right). The mass transportation from the direct wave of the synthetic signal to the reflected wave of the real seismic signal (within the green box of the upper left subgraph) will cause artifacts in the sensitivity kernel, which arise around the reflected wave path (the blue dashed lines of the upper right subgraph). In the lower subgraphs, we can obtain the satisfactory OT map and sensitivity kernel since only direct waves are picked.} 
	\label{fig:section3_52_2_map_kernel}
\end{figure}

The above difference between the real seismic signal and the synthetic signal is further magnified by the squaring scaling. It leads to unreasonable mass transportation from the direct wave of the synthetic signal to the reflected wave of the real seismic signal (upper left subgraph of Figure \ref{fig:section3_52_2_map_kernel}). Therefore, there will be artifacts in the sensitivity kernel $K_i(\bx)$, as we illustrate in the upper right subgraph of Figure \ref{fig:section3_52_2_map_kernel}. On the other hand, if we only consider the direct waves for inversion, the above-mentioned difficulties will be easily solved, as we illustrate in the lower subgraphs of Figure \ref{fig:section3_52_2_map_kernel}. 

\begin{remark}
In fact, the reflected wave signals are also important to constrain the underground velocity structures \cite{HuYaTo:16}. The reflection phases can also be similarly picked, processed, and used for inversion by our approach. However, the utilization of the reflected wave is not trivial, and more technical details are required in practice \cite{BrOpVi:15, XuWaChLaZh:12, ZhBrOpVi:15}. Thus, we will not discuss the issues of the reflected wave in the following sections.
\end{remark}

\subsection{New normalization method}

As it is well known, the quadratic Wasserstein metric measures the difference between two probability density functions, which is not directly suitable for seismic signals. Thus, some processing procedures, i.e., choosing an appropriate operator $\mathcal{P}$ in \eqref{eqn:ref_dis} are required to convert the seismic signals into probability density functions. Several different approaches, e.g., linear scaling \cite{YaEnSuFr:18}, squaring scaling \cite{ChChWuYa:18}, and exponential scaling \cite{QiJaVaYaEn:17}, have been proposed to address this issue. Among these methods, the squaring scaling maintains convexity very well, and it is worthy of more discussions.

The normalization operator with squaring scaling consists of two ingredients: squaring seismic signal to ensure non-negativity and normalization to guarantee the same mass. A natural approach is

\begin{equation}
	\mathcal{P}_{1}(s(t))=\frac{s^2(t)}{\norm{s^2(t)}},
	\label{equ:normalization1}
\end{equation}
in which

\begin{equation*}
     \norm{s(t)}=\int_{0}^{t_f} s(t) \rd t.
\end{equation*}

Substitute the above formula into equation \eqref{eqn:ref_dis}, the form of the misfit function is given by

\begin{equation*}
\chi=\mcd(s(t),d(t))=W^2_2\left(\frac{s^2(t)}{\norm{s^2(t)}},\frac{d^2(t)}{\norm{d^2(t)}}\right).
\end{equation*}
Here the subscript indices $i$ and $j$ are dropped for simplicity. According to the discussions in Section \ref{sec:Wasserstein}, we need to compute the inverse of the following cumulative distribution function

\begin{equation*}
G(t) = \int_0^{t_f} \frac{d^2(t)}{\norm{d^2(t)}} \rd t.
\end{equation*}
However, $G^{-1}(t)$ is not well defined when the real seismic signal $d(t)=0$ in certain interval. Correspondingly, there will be difficulties in the computation of the misfit function.

\begin{figure}[H]
	\centering  
	\includegraphics[width=1\textwidth]{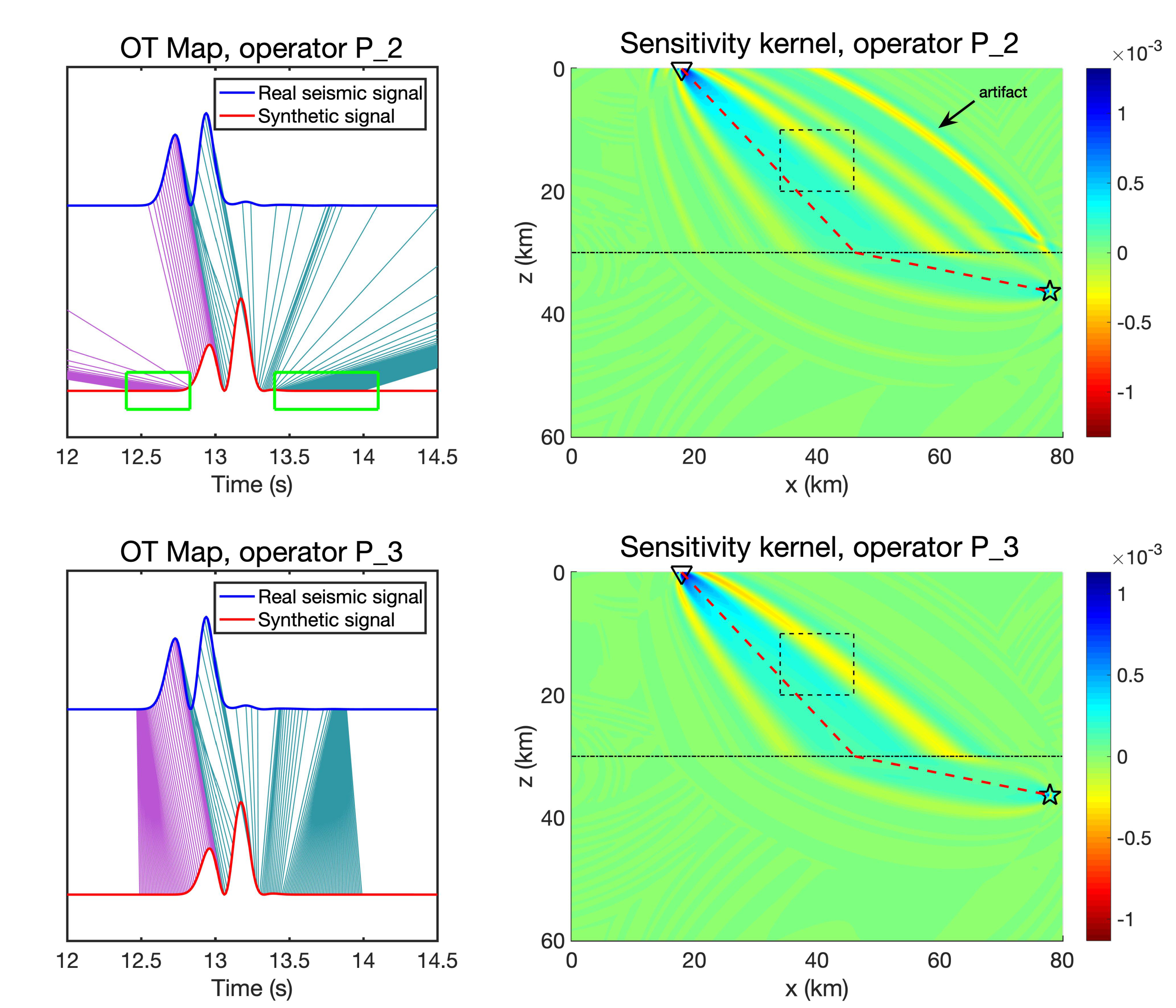}
	\caption{Illustration of the Optimal Transport map between the real seismic signal and synthetic signal (left) and the sensitivity kernel (right). In the upper subgraphs, the newly created mass by the operator $\mathcal{P}_{2}$ could not be balanced, which leads to unreasonable mass transportation (upper left) and artifacts in the sensitivity kernel (upper right). In the lower subgraphs, we can obtain the satisfactory OT map and sensitivity kernel since a new operator $\mathcal{P}_{3}$ is used. }
	\label{fig:section3_79_6_map_kernel}
\end{figure}

In order to avoid the above-mentioned problem, we can make a slight upward shift on the squared signal before the normalization, i.e., 

\begin{equation}
	\mathcal{P}_{2}(s(t))=\frac{s^2(t)+\varepsilon}{\norm{s^2(t)+\varepsilon }}.
	\label{equ:normalization2}
\end{equation}
Here $\varepsilon>0$ is a small parameter. However, the misfit function in \eqref{eqn:ref_dis} with this normalization operator

\begin{equation*}
\chi=\mcd(s(t),d(t))=W^2_2\left( \frac{s^2(t)+\varepsilon}{\norm{s^2(t)+\varepsilon}},\frac{d^2(t)+\varepsilon}{\norm{d^2(t)+\varepsilon}} \right)
\end{equation*}
still leads to unreasonable mass transportation (green box in the upper left subgraph of Figure  \ref{fig:section3_79_6_map_kernel}) since the additional mass does not equal

\begin{equation*}
\frac{\varepsilon}{\norm{s^2(t)+\varepsilon}} \neq \frac{\varepsilon}{\norm{d^2(t)+\varepsilon}}.
\end{equation*}
This again leads to artifacts in the sensitivity kernel $K_i(\bx)$ (upper right subgraph of Figure
\ref{fig:section3_79_6_map_kernel}).

With a simple trick, we can solve the problem of unequal additional masses by modifying the normalization operator as

\begin{equation}
	\mathcal{P}_{3}(s(t))=\frac{\frac{s^2(t)}{\norm{s^2(t)}}+\varepsilon}{1+t_f\varepsilon}.
	\label{equ:newnormalization}
\end{equation} 
We can clearly see that regardless of the values of $s(t)$ and $d(t)$, the additional mass is $\frac{\varepsilon}{1+t_f \varepsilon}$. As a result, we can avoid all the mentioned troubles. Both the OT map and the sensitivity kernel are satisfactory, as we illustrate in the lower subgraphs of Figure \ref{fig:section3_79_6_map_kernel}.

\begin{remark}
	In the squaring scaling, a parameter $\varepsilon$ is added to avoid the singularity. It is noted that large $\varepsilon$ could destroy the convexity property. On the other hand, there will still be numerical singularities when $\varepsilon$ is small. In practice, $\varepsilon$ is feasible in a relatively large range, e.g., $10^{-4} \sim 10^{-2}$. In the following numerical experiments, we select $\varepsilon = 10^{-3}$. 
\end{remark}

\section{Numerical Experiments} \label{sec:num}
In this section, we present two numerical experiments to investigate the validity of our inversion method based on the quadratic Wasserstein metric with squaring scaling. We use the finite difference method to solve the acoustic wave equation \cite{Da:86, LiYaWuMa:17, YaEnSuFr:18}. The perfectly matched layer boundary condition \cite{KoTr:03} is applied to absorb the outgoing wave. The delta source function is discretized by piecewise polynomial given in \cite{We:08}

\begin{equation*}
	\delta_h(x)=\left\{\begin{array}{ll}
		\frac{1}{h}\left(1-\frac{5}{4}\abs{\frac{x}{h}}^2-\frac{35}{12}\abs{\frac{x}{h}}^3
			+\frac{21}{4}\abs{\frac{x}{h}}^4-\frac{25}{12}\abs{\frac{x}{h}}^5\right), & \abs{x}\le h, \\
		\frac{1}{h}\left(-4+\frac{75}{4}\abs{\frac{x}{h}}-\frac{245}{8}\abs{\frac{x}{h}}^2+\frac{545}{24}\abs{\frac{x}{h}}^3
			-\frac{63}{8}\abs{\frac{x}{h}}^4+\frac{25}{24}\abs{\frac{x}{h}}^5\right), & h<\abs{x}\le 2h, \\
		\frac{1}{h}\left(18-\frac{153}{4}\abs{\frac{x}{h}}+\frac{255}{8}\abs{\frac{x}{h}}^2-\frac{313}{24}\abs{\frac{x}{h}}^3
			+\frac{21}{8}\abs{\frac{x}{h}}^4-\frac{5}{24}\abs{\frac{x}{h}}^5\right), & 2h<\abs{x}\le 3h, \\
		0, & \abs{x}>3h.
	\end{array}\right. 
\end{equation*}
Here $h$ is related to the mesh size.
 
\subsection{The Two-Layer Model} \label{subsec:num1}
Consider the two-layer model in a bounded domain $\Omega=[0,80\; km]\times[0,60\; km]$, which consists of the crust, the uppermost mantle, and the Moho discontinuity at a depth of $30\; km$, see Figure \ref{fig:section3_velocity_52_2} for illustration. The real seismic velocity model includes a $+15\%$ high-velocity anomaly in the crust, given by

\begin{equation*}
	c_{T}(x,z)=\left\{\begin{array}{ll}
	        6.67\ km/s, & (x,z) \in [35\ km,45\ km] \times [10\ km,20\ km], \\ 
		8.1\ km/s, & z>30\ km, \\
		5.8\ km/s, & others.
	\end{array}\right.
\end{equation*}
Our goal is to perform the seismic velocity inversion to detect this high-velocity anomaly. Correspondingly, the initial velocity model without high-velocity anomaly is as follows

\begin{equation*}
	c_{0}(x,z)=\left\{\begin{array}{ll} 
		5.8\ km/s, & z \leq 30\ km, \\
		8.1\ km/s, & z>30\ km. 
	\end{array}\right.
\end{equation*}
The computational time interval is $[0\; s,21\; s]$. The inversion grid step is $2\; km$ and the number of degrees of freedom amounts to $1200$. The space and time steps in the forward simulation are $0.2\ km$ and $0.01\ s$, respectively. The dominant frequency of the earthquakes in \eqref{eqn:ricker} is $f_0=2\; Hz$. We randomly choose $25$ receiver stations deployed on the surface and $80$ earthquakes distributed in the study region.

We then perform the seismic velocity inversion by using the quadratic Wasserstein metric with squaring scaling. As a comparison, the inversion is also performed with the traditional $L^2$ metric.
 
To quantitatively compare the results of different methods, we also compute the relative model error 

\begin{equation*}
RME=\frac{\int_{\Omega} |c_k(\bx)-c_T(\bx)|^2 \rd \bx}{\int_{\Omega} |c_0(\bx)-c_T(\bx)|^2 \rd \bx},
\end{equation*} 
and the relative misfit function 

\begin{equation*}
RMF=\frac{\Xi(c_k(\bx))}{\Xi(c_0(\bx))},
\end{equation*} 
where $c_k(\bx)$ indicates the velocity model in the $k$-th iteration.

 \begin{figure}[H]
	\centering
	\includegraphics[width=1\textwidth]{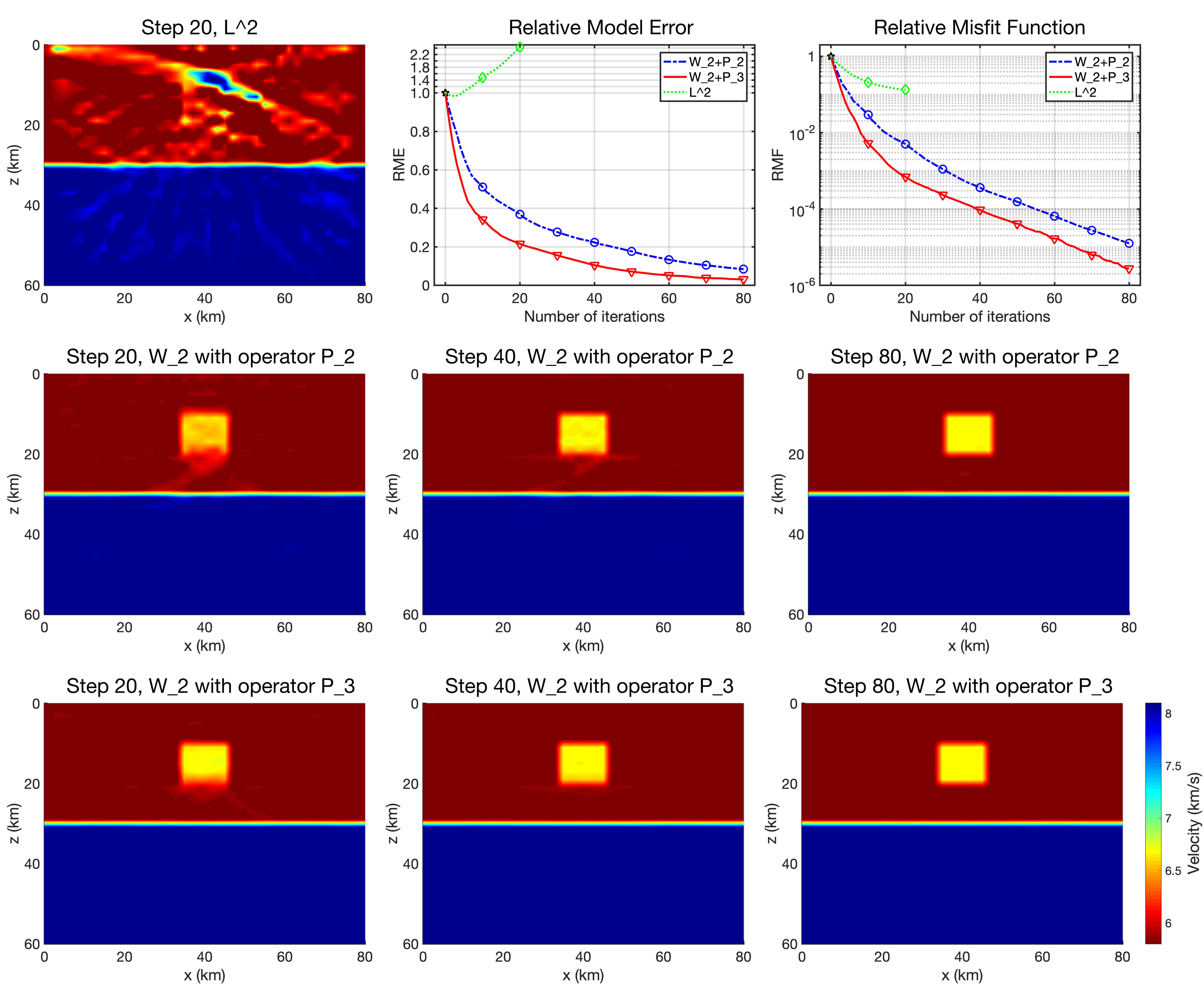}
	\caption{The inversion results of the two-layer model. Upper subgraphs: the result for $L^2$ metric after 20 steps (upper left); the convergent trajectories of the relative model error (upper middle); the convergent trajectories of the relative misfit function (upper right). In the middle and the lower subgraphs, we present the results for the $W_2$ metric with the operators $\mathcal{P}_{2}$ and $\mathcal{P}_{3}$, respectively. From left to right, the inversion iteration steps are $20$, $40$, and $80$. All the results are shown in the same color bar.}  
	\label{fig:exam1_res_tl}
\end{figure}

In Figure \ref{fig:exam1_res_tl}, we present the inversion results of $L^2$ metric and $W_2$ metric. Obviously, the $L^2$-based inversion could not capture the $+15\%$ high-velocity anomaly (upper left subgraph of Figure \ref{fig:exam1_res_tl}). Although the misfit function decreases in the iteration (upper middle subgraph of Figure \ref{fig:exam1_res_tl}), the model error increases (upper right subgraph of Figure \ref{fig:exam1_res_tl}).

In Figure \ref{fig:exam1_res_tl} and Table \ref{tab:exam1_res_tl}, we also compare the inversion results of the quadratic Wasserstein metric with different operators $\mathcal{P}_{2}$ and $\mathcal{P}_{3}$. From the convergent trajectories (upper middle and upper right subgraphs of Figure \ref{fig:exam1_res_tl}), we can see the relative model error and the relative misfit function of the operator $\mathcal{P}_{3}$ both have a faster descent rate than those of the operator $\mathcal{P}_{2}$. Quantitatively, we can see from Table \ref{tab:exam1_res_tl} that the operator $\mathcal{P}_{3}$ only needs half of the iteration steps of the operator $\mathcal{P}_{2}$ to achieve almost the same relative model error and relative misfit function. This significantly saves the expensive computational cost of the seismic velocity inversion problem. Finally, it can be seen from the middle and lower subgraphs of Figure \ref{fig:exam1_res_tl}, the velocity inversion results of the operator $\mathcal{P}_{3}$ are significantly better than those of the operator $\mathcal{P}_{2}$ under the same iteration steps. The above discussions show that our approach has higher efficiency and better inversion results.

\linespread{1.5}
\begin{table}
	\setlength{\belowcaptionskip}{0.3cm}
	\centering
	\caption{The two-layer model. Relative Model Error and Relative Misfit Function of $W_2$ with the operators $\mathcal{P}_{2}$ and $\mathcal{P}_{3}$ in $20$, $40$ and $80$ iteration steps, respectively.}
	\label{tab:exam1_res_tl}
	\begin{tabular}{ccccc}
	\toprule
	\multirow{2}{*}{Iteration Steps} & \multicolumn{2}{c}{Relative Model Error} & \multicolumn{2}{c}{Relative Misfit Function} \\ 
	\cline{2-5}
	         & $W_2$ with $P_2$ & $W_2$ with $P_3$ & $W_2$ with $P_2$ & $W_2$ with $P_3$ \\
         \midrule	
	 $20$ & $3.69 \times 10^{-1}$ & $2.15 \times 10^{-1}$ & $4.99 \times 10^{-3}$ & $6.90 \times 10^{-4}$ \\ 

	 $40$ & $2.23 \times 10^{-1}$ & $1.04 \times 10^{-1}$ & $3.61 \times 10^{-4}$ & $9.41 \times 10^{-5}$ \\ 
	 
	  $80$ & $ 8.35 \times 10^{-2}$ & $3.04 \times 10^{-2}$ & $1.25 \times 10^{-5}$ & $2.75 \times 10^{-6}$ \\ 
	   \bottomrule
	\end{tabular}
\end{table}

\subsection{The Crustal Root Model} \label{subsec:num2}

\begin{figure}[H]
	\centering
	\includegraphics[width=0.82\textwidth]{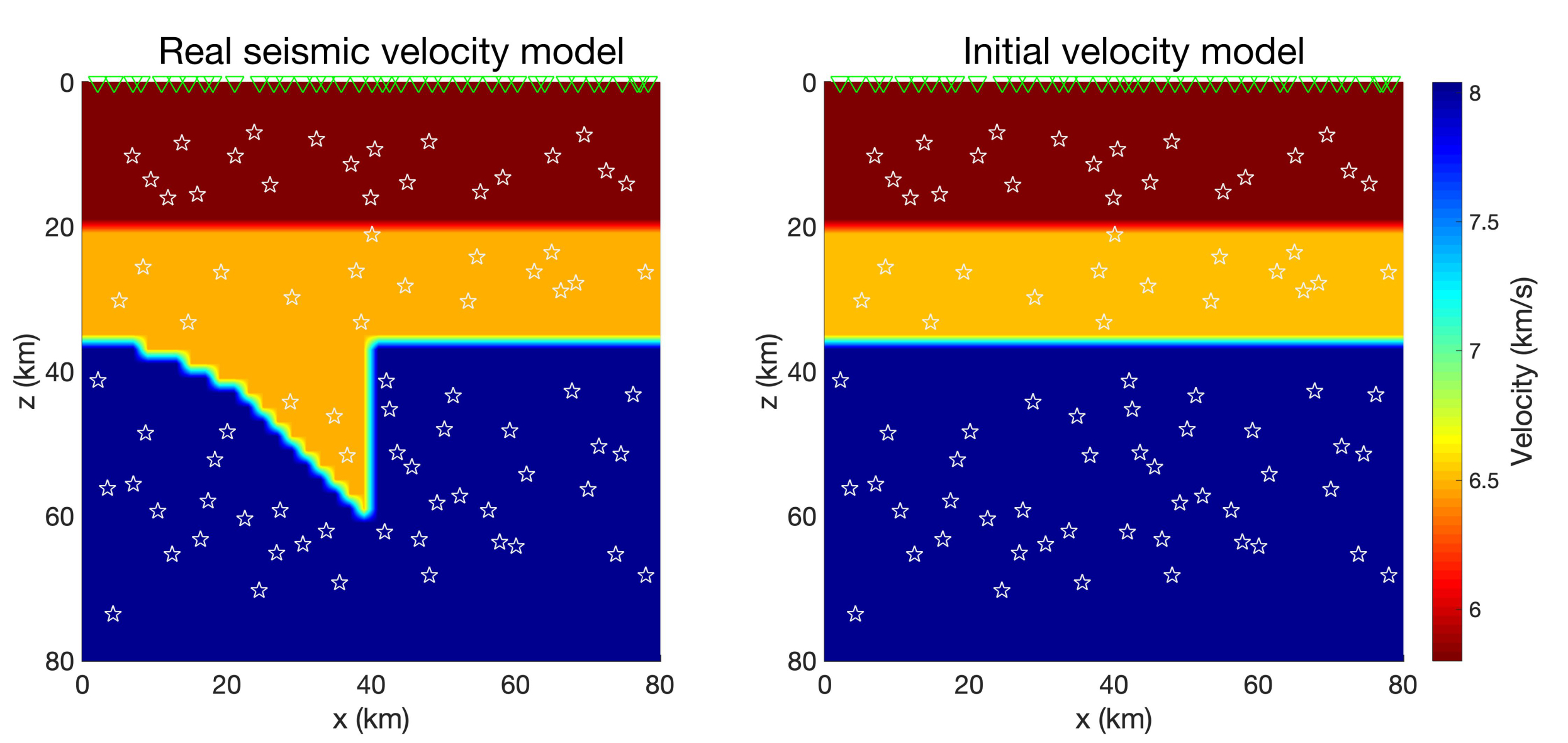}
	\caption{Illustration of the crustal root model. Left: the real seismic velocity model. Right: the initial velocity model. The green inverted triangles and the white stars indicate the receiver stations and the earthquakes, respectively.}
	\label{fig:exam2_back_velocity_cr}
\end{figure}

Let us consider the crustal root model, a kind of subsurface structure usually found along the orogen. This model consists of the two-layered crust divided by the Conrad discontinuity. A dipping and discontinuous Moho interface separates the crust and the mantle. The depiction of these tectonic features helps us better understand the forming of the old mountains.  In mathematics, we consider this three-layer model in the bounded domain $\Omega=[0,80\; km]\times[0,80\; km]$. Three layers are divided by the Conrad discontinuity at $20\; km$ depth and the Moho discontinuity whose location $(x, L(x))$ is formulated with a quadratic function is given by

\begin{equation*}
L(x)=\left\{\begin{array}{ll}
            36+\frac{25}{1600}x^2\ km , & 0\ km \le x \le 40\ km, \\
            36\ km , & 40\ km < x \le 80\ km.
            \end{array}\right.
\end{equation*}
The seismic wave speed at each layer refers to the AK135 model \cite{KeEnBu:95}, generating the real seismic velocity model (Figure \ref{fig:exam2_back_velocity_cr}, left)

\begin{equation*}
c_{T}(x,z)=\left\{\begin{array}{ll}
	        5.8\ km/s, & z\le 20\ km, \\ 
		6.5\ km/s, & 20\ km < z\le L(x), \\
		8.04\ km/s, & others.
	\end{array}\right.
\end{equation*}
Our goal is to perform the seismic velocity inversion to detect this crustal root. Correspondingly, the initial velocity model (Figure \ref{fig:exam2_back_velocity_cr}, right) without crustal root anomaly is as follows

\begin{equation*}
c_{0}(x,z)=\left\{\begin{array}{ll}
	        5.8\ km/s, & z\le 20\ km, \\ 
		6.5\ km/s, & 20\ km < z\le 36\ km \\
		8.04\ km/s, & others.
	\end{array}\right.
\end{equation*}
The computational time interval is $[0\; s,21\; s]$. The inversion grid step is $2\; km$ and the number of degrees of freedom amounts to $1600$. The space and time steps in the forward simulation are $0.2\ km$ and $0.01\ s$, respectively. The dominant frequency of the earthquakes in \eqref{eqn:ricker} is $f_0=2\; Hz$. We randomly choose $40$ receiver stations deployed on the surface and $80$ earthquakes distributed in the study region.

\begin{figure}
	\centering
	\includegraphics[width=1\textwidth]{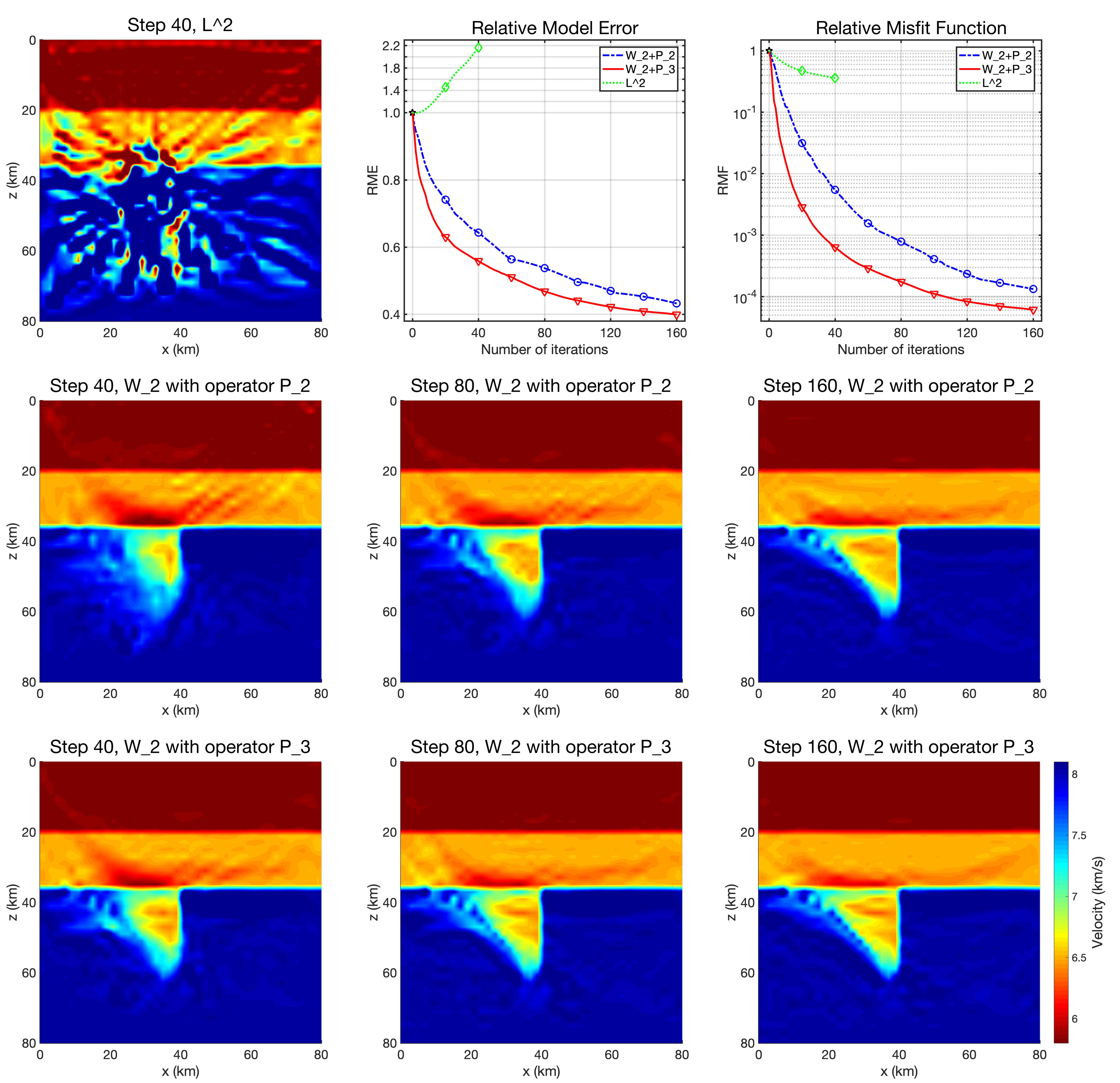}
	\caption{The inversion results of the crustal root model. Upper subgraphs: the result for $L^2$ metric after 40 steps (upper left); the convergent trajectories of the relative model error (upper middle); the convergent trajectories of the relative misfit function (upper right). In the middle and the lower subgraphs, we present the results for the $W_2$ metric with the operators $\mathcal{P}_{2}$ and $\mathcal{P}_{3}$, respectively. From left to right, the inversion iteration steps are $40$, $80$, and $160$. All the results are shown in the same color bar.}  
	\label{fig:exam2_res_cr}
\end{figure}

Similar to subsection \ref{subsec:num1}, we present the inversion results of $L^2$ metric and $W_2$ metric with the operators $\mathcal{P}_{2}$ and $\mathcal{P}_{3}$ in Figure \ref{fig:exam2_res_cr}. Obviously, the $L^2$-based inversion could not capture the crustal root structure. The relative model error and the relative misfit function with respect to different normalization operators are given in Table \ref{tab:exam2_res_cr}. Correspondingly, the convergent trajectories are output in the upper middle and upper right subgraphs of Figure \ref{fig:exam2_res_cr}. In the middle and lower subgraphs of Figure \ref{fig:exam2_res_cr}, the inversion results are also presented. From which, we can draw the same conclusions as those in subsection \ref{subsec:num1}.

\linespread{1.5}
\begin{table}
	\setlength{\belowcaptionskip}{0.3cm}
	\centering
	\caption{The crustal root model. Relative Model Error and Relative Misfit Function of $W_2$ with the operators $\mathcal{P}_{2}$ and $\mathcal{P}_{3}$ in $40$, $80$ and $160$ iteration steps, respectively.}
	\label{tab:exam2_res_cr}
	\begin{tabular}{cccccc}
	\toprule
	\multirow{2}{*}{Iteration Steps} & \multicolumn{2}{c}{Relative Model Error} & \multicolumn{2}{c}{Relative Misfit Function} \\ 
	\cline{2-5}
	         & $W_2$ with $P_2$ & $W_2$ with $P_3$ & $W_2$ with $P_2$ & $W_2$ with $P_3$ \\
       \midrule	
	 $40$ & $6.43 \times 10^{-1}$ & $5.59 \times 10^{-1}$ & $5.47 \times 10^{-3}$ & $6.35 \times 10^{-4}$ \\ 

	 $80$ & $5.37 \times 10^{-1}$ & $4.68 \times 10^{-1}$ & $7.83 \times 10^{-4}$ & $1.74 \times 10^{-4}$ \\ 

	  $160$ & $ 4.32 \times 10^{-1}$ & $3.99 \times 10^{-1}$ & $1.33 \times 10^{-4}$ & $6.11 \times 10^{-5}$ \\ 
	  \bottomrule
	\end{tabular}
\end{table}

\section{Conclusion} \label{sec:con}

What we have seen from the above is the solution to the problem that the seismic velocity inversion based on squaring scaling and the quadratic Wasserstein metric is difficult, as mentioned in [Commun. Inf. Syst., 2019, 19:95-145] and [Meth. Appl. Anal., 2019, 2:133-148]. We can not only solve the seismic velocity inversion with a large number of degrees of freedom. By introducing a better normalization operator, the convergence efficiency is significantly improved. We would like to combine the above techniques with the double-difference traveltime adjoint tomography \cite{ChChWuYaTo:21}, which has significant advantages in real seismic data. This may result in a more robust and reliable seismic velocity inversion method. We are currently investigating this interesting topic and hope to report this in an independent publication.

\section*{Acknowledgments}
This work was supported by National Natural Science Foundation of China (Grant No. 11871297) and Tsinghua University Initiative Scientific Research Program.



\begin{thebibliography}{}
 \bibitem{AkRi:80}
	K. Aki and P.G. Richards, \textit{Quantitative Seismology: Theory and Methods volume II}, W.H. Freeman \& Co (Sd), 1980.

      \bibitem{BrOpVi:15}
      R. Brossier, S. Operto and J. Virieux, Velocity model building from seismic reflection data by full-waveform inversion, \textit{Geophysical Prospecting}, \textbf{63}(2), 354-367, 2015.
      
	\bibitem{ChBa:06}
	S.-J. Chang and C.-E. Baag, Crustal Structure in Southern Korea from Joint Analysis of Regional Broadband Waveforms and Travel Times, \textit{Bulletin of the Seismological Society of America}, \textbf{96}(3), 856–870, 2006.

    \bibitem{ChChWuYaTo:21}
     J. Chen, G.X. Chen, H. Wu, J.Y. Yao and P. Tong, Adjoint tomography of NE Japan revealed by common-source double-difference traveltime data, \textit{preprint}.

       \bibitem{ChChWuYa:18}
	J. Chen, Y.F. Chen, H. Wu and D.H. Yang, The quadratic Wasserstein metric for Earthquake Location, \textit{J. Comput. Phys.}, \textbf{373}, 188-209, 2018.

	\bibitem{ChHiLe:17}
	Y. Chen, J. Hill, W. Lei, M. Lefebvre, J. Tromp, E. Bozdag and D. Komatitsch, Automated time-window selection based on machine learning for full-waveform inversion, \textit{SEG Technical Program Expanded Abstracts}, 1604-1609, 2017.
				  
	\bibitem{ChNiPi:14}
	R. Chu, S. Ni, A. Pitarka and D.V. Helmberger,  Inversion of Source Parameters for Moderate Earthquakes Using Short-Period Teleseismic P Waves, \textit{Pure and Applied  Geophysics}, \textbf{171}(7), 1329–1341, 2014.
		
	\bibitem{Da:86}
	M.A. Dablain, The application of high-order differencing to the scalar wave equation, \textit{Geophysics}, \textbf{51}(1), 54-66, 1986.	

       \bibitem{DuYa:20}		
        M. Dunlop and Y.N. Yang, New likelihood functions and level-set prior for Bayesian full-waveform inversion, In \textit{SEG Technical Program Expanded Abstracts 2020}, pages 825-829. Society of Exploration Geophysicists, 2020.

    \bibitem{EnFr:14}
        B. Engquist and B.D. Froese, Application of the Wasserstein metric to seismic signals, \textit{Commun. Math. Sci.}, \textbf{12}(5), 979-988, 2014.

    \bibitem{EnFrYa:16}
        B. Engquist, B.D. Froese and Y.N. Yang, Optimal transport for seismic full waveform inversion, \textit{Commun. Math. Sci.}, \textbf{14}(8), 2309-2330, 2016.
      
      \bibitem{EnYa:19a}
       B. Engquist and Y.N. Yang, Seismic imaging and optimal transport, \textit{Communications in Information and Systems}, \textbf{19}(2), 95-145, 2019.  
        
        \bibitem{EnYa:19b}
        B. Engquist and Y.N. Yang, Seismic inversion and the data normalization for optimal transport, \textit{Methods and Applications of Analysis}, \textbf{26}(2), 133-148, 2019.
 
        
     \bibitem{EnYa:21}
       B. Engquist and Y.N. Yang, Optimal Transport Based Seismic Inversion:Beyond Cycle Skipping, \textit{Comm. Pure Appl. Math.}. https://doi.org/10.1002/cpa.21990, 2021.

    \bibitem{HuDaNo:01}
     S.-H. Hung, F.A.Dahlen and G. Nolet, Wavefront healing: a banana–doughnut perspective, \textit{Geophys. J. Int.}, \textbf{146}(2), 289-312, 2001.
     
     \bibitem{HuYaTo:16}
	X. Huang, D. Yang, P. Tong, J. Badal and Q. Liu, Wave equation-based reflection tomography of the 1992 Landers earthquake area, \textit{Geophysical Research Letters}, \textbf{43}(5), 1884-1892, 2016.

	\bibitem{KoTr:03}
	D. Komatitsch and J. Tromp, A perfectly matched layer absorbing boundary condition for the second-order seismic wave equation, \textit{Geophys. J. Int.}, \textbf{154}(1), 146-153, 2003.				 
	\bibitem{KeEnBu:95}
        B.L.N. Kennett, E.R. Engdahl and R. Buland, Constraints on seismic velocities in the Earth from traveltimes, \textit{Geophys. J. Int.}, \textbf{122}(1), 108-124, 1995.

	\bibitem{LiYaWuMa:17}
	J.S. Li, D.H. Yang, H. Wu and X. Ma, A low-dispersive method using the high-order stereo-modelling operator for solving 2-D wave equations, \textit{Geophys. J. Int.}, \textbf{210}(3), 1938-1964, 2017.
				
	\bibitem{Ma:15}
	 R. Madariaga, Seismic Source Theory, in \textit{Treatise on Geophysics (Second Edition)}, pp. 51-71, ed. Gerald, S., Elsevier B.V., 2015.

	\bibitem{MaTaChChTr:09}
	 A. Maggi, C. Tape, M. Chen, D. Chao and J. Tromp, An automated time-window selection algorithm for seismic tomography, \textit{Geophys. J. Int.}, \textbf{178}(1), 257–281, 2009.

       \bibitem{MeBrViOp:13} 
          L. Métivier, R. Brossier, J. Virieux and S. Operto,  Full waveform inversion and the truncated Newton method, \textit{SIAM Journal on Scientific Computing}, \textbf{35}(2), B401-B437, 2013.
      
       \bibitem{MeBrMeOuVi:16a}
        L. Metivier, R. Brossier, Q. Merigot, E. Oudet and J. Virieux, Measuring the misfit between seismograms using an optimal transport distance: application to full waveform inversion, \textit{Geophys. J. Int.}, \textbf{205}(1), 345–377, 2016. 

       \bibitem{MeBrMeOuVi:16b}
        L. Metivier, R. Brossier, Q. Merigot, E. Oudet and J. Virieux, An optimal transport approach for seismic tomography: application to 3D full waveform inversion, \textit{Inverse Problems}, \textbf{32}(11), 115008, 2016. 
     
       \bibitem{Pl:06}
        R.-E. Plessix, A review of the adjoint-state method for computing the gradient of a functional with geophysical applications, \textit{Geophys. J. Int.}, \textbf{167}(2), 495–503, 2006.

       \bibitem{QiJaVaYaEn:17}		
        L. Qiu, J. Ramos-Mart\'inez, A. Valenciano, Y. Yang and B. Engquist, Full-waveform inversion with an exponentially encoded optimal-transport norm, In \textit{SEG Technical Program Expanded Abstracts 2017}, pages 1286–1290. Society of Exploration Geophysicists, 2017.
        
        \bibitem{RaSaHa:10}
        N. Rawlinson, M. Sambridge and J. Hauser, Multipathing, reciprocal traveltime fields and raylets, \textit{Geophys. J. Int.}, \textbf{181}(2), 1077-1092, 2010.
	
        \bibitem{TaLiuMaTr:10}
        C. Tape, Q.Y. Liu, A. Maggi and J. Tromp, Seismic tomography of the southern California crust based on spectral-element and adjoint methods, \textit{Geophys. J. Int.}, \textbf{180}(1), 433-462, 2010.
	
       \bibitem{Vi:03}
        C. Villani, \textit{Topics in Optimal transportation}, Graduate Studies in Mathematics, American Mathematical Society, 2003.

       \bibitem{Vi:08}
        C. Villani, \textit{Optimal Transport: Old and New}, Springer Science $\&$ Business Media, 2008.	 	
	\bibitem{ViOp:09}
        J. Virieux and S. Operto,  An overview of full-waveform inversion in exploration geophysics, \textit{Geophysics}, \textbf{74}(6), WCC1-WCC26, 2009. 
 	
	\bibitem{WaYaJiWu:19} 
        J. Wang, D.H. Yang, H. Jing and H. Wu, Full waveform inversion based on the ensemble Kalman filter method using uniform sampling without replacement, \textit{Science Bulletin},  \textbf{64}(5), 321-330, 2019.
	
	\bibitem{WaTk:20}
	S. Wang and H. Tkalčić, Seismic event coda-correlation: Toward global coda-correlation tomography, \textit{Journal of Geophysical Research: Solid Earth}, \textbf{125}, e2019JB018848, 2020.
			
	\bibitem{We:08}
   		X. Wen, High Order Numerical Quadratures to One Dimensional Delta Function Integrals, \textit{SIAM J. Sci. Comput.}, \textbf{30}(4), 1825-1846, 2008.


    \bibitem{XuWaChLaZh:12}
    S. Xu, D.L. Wang, F. Chen, G. Lambar\'e and Y. Zhang, Inversion on Reflected Seismic Wave, \textit{SEG Technical Program Expanded Abstracts}, 1-7, 2012.
 		
    \bibitem{YaEnSuFr:18}
        Y.N. Yang, B. Engquist, J.Z. Sun and B.D. Froese, Application of Optimal Transport and the Quadratic Wasserstein Metric to Full-Waveform Inversion, \textit{Geophysics}, \textbf{83}(1), R43-R62, 2018.

    \bibitem{YaEn:18}
        Y.N. Yang and B. Engquist, Analysis of optimal transport and related misfit functions in full-waveform inversion, \textit{Geophysics}, \textbf{83}(1), A7-A12, 2018.

    \bibitem{ZhChWuYaQi:21}
       D. T. Zhou, J. Chen, H. Wu, D.H. Yang and L.Y. Qiu, The Wasserstein-Fisher-Rao metric for waveform based earthquake location, \textit{J. Comp. Math.}, accepted.	

    \bibitem{ZhBrOpVi:15}
     W. Zhou, R. Brossier, S. Operto and J. Virieux, Full waveform inversion of diving \& reflected waves for velocity model building with impedance inversion based on scale separation, \textit{Geophys. J. Int.}, \textbf{202}(3), 1535-1554, 2015.

\end{thebibliography}
\end{document}